\documentclass [11pt]{article}
\usepackage{mathrsfs}
\usepackage{amsthm,amsmath,amsfonts,amssymb,graphicx,graphics,psfrag}

\newtheorem{thm}{Theorem}

\newtheorem{lem}{Lemma}

\date{}

\begin{document}

\title{Uniform Star-factors of Graphs with Girth Three
\thanks{This work is supported by RFDP of Higher Education of China and Natural Sciences
and Engineering Research Council of Canada.}}

\author{ Yunjian Wu$^{1}$ and Qinglin Yu$^{1,2}$
\thanks{Corresponding author: yu@nankai.edu.cn}
\\ {\small $^1$ Center for Combinatorics, LPMC}
\\ {\small \small Nankai University, Tianjin, 300071, China}
\\ {\small  \small $^2$ Department of Mathematics and Statistics}
\\ {\small  Thompson Rivers University, Kamloops, BC, Canada}
}

\maketitle

\begin{abstract}
 A {\it star-factor} of a graph $G$ is a spanning subgraph of $G$
such that each component of which is a star. Recently, Hartnell and
Rall studied a family $\mathscr{U}$ of graphs satisfying the
property that every star-factor of a member graph has the same
number of edges. They determined the family $\mathscr{U}$ when the
girth is at least five. In this paper, we investigate the family of
graphs with girth three and determine all members of this family.
\begin{flushleft}
{\em Key words:} star-factor, uniform star-factor, girth, edge-weighting \\
\end{flushleft}
\end{abstract}

\vskip 1mm

\vskip 1mm

\section{Introduction}
Throughout this paper, all graphs considered are simple. We refer
the reader to \cite{bb} for standard graph theoretic terms not
defined in this paper.

Let $G=(V,E)$ be a graph with vertex set $V(G)$ and edge set $E(G)$.
If $G$ is not a forest, the length of a shortest cycle in $G$ is
called the {\it girth} of $G$. We say that a forest has an infinite
girth.  We shall often construct new graphs from old ones by
deleting some vertices or edges. If $W\subset V(G)$, then
$G-W=G[V-W]$ is the subgraph of $G$ obtained by deleting the
vertices in $W$ and all edges incident with them. Similarly, if
$E'\subset E(G)$, then $G-E'=(V(G), E(G)-E')$. We denote the degree
of a vertex $x$ in $G$ by $d_{G}(x)$, and the set of vertices
adjacent to $x$ in $G$ by $N_{G}(x)$. A {\it leaf} is a vertex of
degree one and a {\it stem} is a vertex which has at least one leaf
as its neighbor. A {\it star} is a tree isomorphic to $K_{1,n}$ for
some $n \geq 1$, and the vertex of degree $n$ is called the {\it
center} of the star. A {\it star-factor} of a graph $G$ is a
spanning subgraph of $G$ such that each component of which is a
star. Clearly a graph with isolated vertices has no star-factors. It
is not hard to see that every graph without isolated vertices admits
a star-factor. If one limits the size of the star used, the
existence of such a star-factor is non-trivial. In \cite{amahashi},
Amahashi and Kano presented a criterion for the existence of a
star-factor, i.e.,  \{$K_{1, 1}, \cdots , K_{1, n}$\}-factor. Yu
\cite{yu} obtained an upper bound on the number of edges in a graph
with unique star-factor.

An {\it edge-weighting} of a graph $G$ is a function $w:
E(G)\longrightarrow \mathbb{N}^+$, where $\mathbb{N}^+$ is the set
of positive integers. For a subgraph $H$, the {\it weight} of $H$
under $w$ is the sum of all the weight values for edges belonging to
$H$, i.e., $w(H)=\Sigma_{e\in E(H)}w(e)$. Motivated by the minimum
cost spanning tree and the optimal assignment problems, Hartnell and
Rall posed an interesting general question: for a given graph, does
there exist an edge-weighting function $w$ such that a certain type
of spanning subgraphs always has the same weights.  In particular,
they investigated the following narrow version of the problem in
which the spanning subgraph is a star-factor.\\

\noindent\textbf{Star-Weighting Problem {{\rm{(Hartnell and Rall
\cite{hartnell}}}}}): For a given graph $G=(V, E)$, is there an
edge-weighting $w$ of $G$ such that every star-factor of $G$ has the
same weights under $w$?

\vspace{4mm}

  To start the investigation, one may consider that the special case that
$w$ is a constant function, i.e., all edges in $G$ are assigned with
the same weights. In this case, every star-factor of $G$ has the
same weights if and only if all star-factors have the same number of
edges. For simplicity, we assume that all edges are assigned with
weight one.

We denote by $\mathscr{U}$ the family of all graphs $G$ such that if
$S_{1}$ and $S_{2}$ are any two star-factors of $G$, then $S_{1}$
and $S_{2}$ have the same number of edges. Clearly, $S_{1}$ and
$S_{2}$ have the same number of edges is equivalent to that they
have the same number of components. Hartnell and Rall classified the
family $\mathscr{U}$ when graphs in $\mathscr{U}$ have girth at
least five and minimum degree at least two.

\begin{thm}
{\rm{(Hartnell and Rall [3])}} Let $G$ be a connected graph of girth
at least five and minimum degree at least two. Then all star-factors
of $G$ have the same weights if and only if $G$ is a $5$-cycle or
$7$-cycle.
\end{thm}

In this paper, we investigate the family $\mathscr{U}$ with girth
three and minimum degree at least two, and we are able to determine
this family completely. The main theorem is as follows.

\begin{thm}
Let $G$ be a connected graph of girth three and minimum degree at
least two. Then all star-factors of $G$ have the same weights if and
only if $G$ is one of the five graphs shown in Figure 1.
\end{thm}

\begin{figure}[h,t]
\begin{center}
\includegraphics[width=7cm]{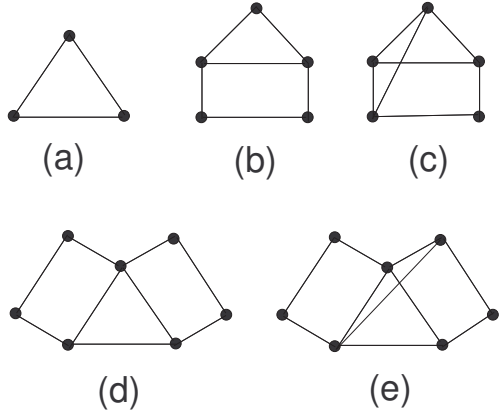}
\end{center}
\vspace{-10mm}
 \caption{Graphs in $\mathscr{U}$ with girth three and
minimum degree at least two}
\end{figure}

\section{Proof of Theorem 2}

Note that if $H$ is a spanning subgraph of $G$, then any star-factor
of $H$ is also a star-factor of $G$.  The following lemma will be
used frequently in reducing the problem of determining membership in
$\mathscr{U}$ to its spanning subgraphs.

\begin{lem}
{\rm{(Hartnell and Rall [3])}} Let $F$ be a subset of $E(G)$ such
that $G-F$ has no isolated vertices.  If $G-F$ is not in
$\mathscr{U}$, then $G$ is not in $\mathscr{U}$.
\end{lem}

    The above lemma implies that if $G$ is in $\mathscr{U}$, then so is
$G-F$.

    The idea to show that a graph does not belong to $\mathscr{U}$ is
to decompose $G$ into several components without isolated vertices
and then simply find one of them not belonging to $\mathscr{U}$.
    For the proof of Theorem $2$, we shall also use the following two lemmas.

\begin{lem}
Let $G$ be a graph with a triangle such that two of its vertices are
of degree two and the third is a stem, then $G$ does not belong to
$\mathscr{U}$.
\end{lem}

\noindent {\bf Proof.} Let $C_{3}=v_{1}v_{2}v_{3}$ be a triangle of
$G$, where $v_{1}$ is a stem adjacent to a leaf $u$. Let $S$ be a
star-factor of the graph $G-\{v_{2},v_{3}\}$. Note that $v_{1}$ is
the center of some star $T$ in $S$. Let $T'$ be the star formed from
$T$ by adding leaves $v_{2}$ and $v_{3}$ adjacent to $v_{1}$, and
let $S'=(S-\{T\}) \cup \{T'\}$, then $S'$ as well as $S\cup
\{v_{2}v_{3}\}$ are star-factors of $G$ having different weights.
Hence $G\notin \mathscr{U}$.\hfill$\Box$

\vspace{2mm}
\begin{lem}
Let $G$ be a graph in $\mathscr{U}$ with a triangle. If exactly
one of the vertices on this triangle has degree at least three,
then all of its neighbors that don't belong to this triangle must
be stems.
\end{lem}

\noindent {\bf Proof.} Suppose $G$ is a graph satisfying the
hypothesis. Let $v$ be a vertex on the triangle of degree at least
three and assume $v$ has a neighbor $x$ not on the triangle such
that $x$ is not a stem. By Lemma $2$, $x$ is not a leaf. Let $F$ be
the set of edges not including $vx$ that are incident with $x$. The
graph $G-F$ has no isolated vertices, and the vertex $v$ is a stem
belonging to a triangle of the type that satisfies the hypothesis of
Lemma $2$, then $G-F$ is not in $\mathscr{U}$. Thus $G$ does not
belong to $\mathscr{U}$ by Lemma $1$, a contradiction.\hfill$\Box$\\

    Now we proceed to prove our main result.\\

\noindent {\bf Proof of Theorem 2.}\ \ \ The only star-factor of a
triangle $C_{3}$ has weight two, so $C_{3}\in \mathscr{U}$. Assume
$G$ belongs to $\mathscr{U}$ and has girth three and minimum degree
at least two but $G$ is not a triangle. Then $G$ contains a triangle
$C_{3}$ with at least two vertices of degree at least three by Lemma
$3$.

Let $C_{3}=v_{1}v_{2}v_{3}$. We consider the following two cases.

\vspace{2mm}

{\it Case $1$.} $d_{G}(v_{3})=2$, $d_{G}(v_{1})\geq 3$,
$d_{G}(v_{2})\geq 3$.

\vspace{2mm}

Let $F_{1}$ be the set of edges incident with $v_{1}$ except
$v_{1}v_{2}$ and $v_{1}v_{3}$. Then $G-F_{1}\in \mathscr{U}$ by
Lemma $1$ since no isolated vertices created in $G-F_{1}$, and all
neighbors of $v_{2}$ not in the triangle $C_{3}$ are stems in
$G-F_{1}$ by Lemma $3$. Let $x$ be a neighbor of $v_{2}$, then there
exists a leaf $y$ incident with $x$. By the definition of $G-F_{1}$,
$y$ is adjacent to $v_{1}$ in $G$ and $d_{G}(y)=2$. Let $F_{2}$ be
the set of edges incident with $v_{2}$ except $v_{2}v_{1}$ and
$v_{2}v_{3}$. A similar argument yields that all neighbors of
$v_{1}$ are stems in $G-F_{2}$ and so $y$ is a stem in $G-F_{2}$.
However $d_{G}(y)=2$, so $x$ is the only leaf of $y$ in $G-F_{2}$.
It follows that $d_{G}(x)=2$ and $x$, $y$, $v_{1}$ and $v_{2}$ form
a quadrangle in $G$. From the above discussion, we see that all
neighbors of $v_{1}$ and $v_{2}$ except $v_1$ and $v_2$ are of
degree two in $G$. Hence $G$ is isomorphic to the graph shown in
Figure $2(a)$ (dashed line indicates a possible edge).

\vspace{3mm}

\begin{figure}[h,t]
\begin{center}
\includegraphics[width=8cm]{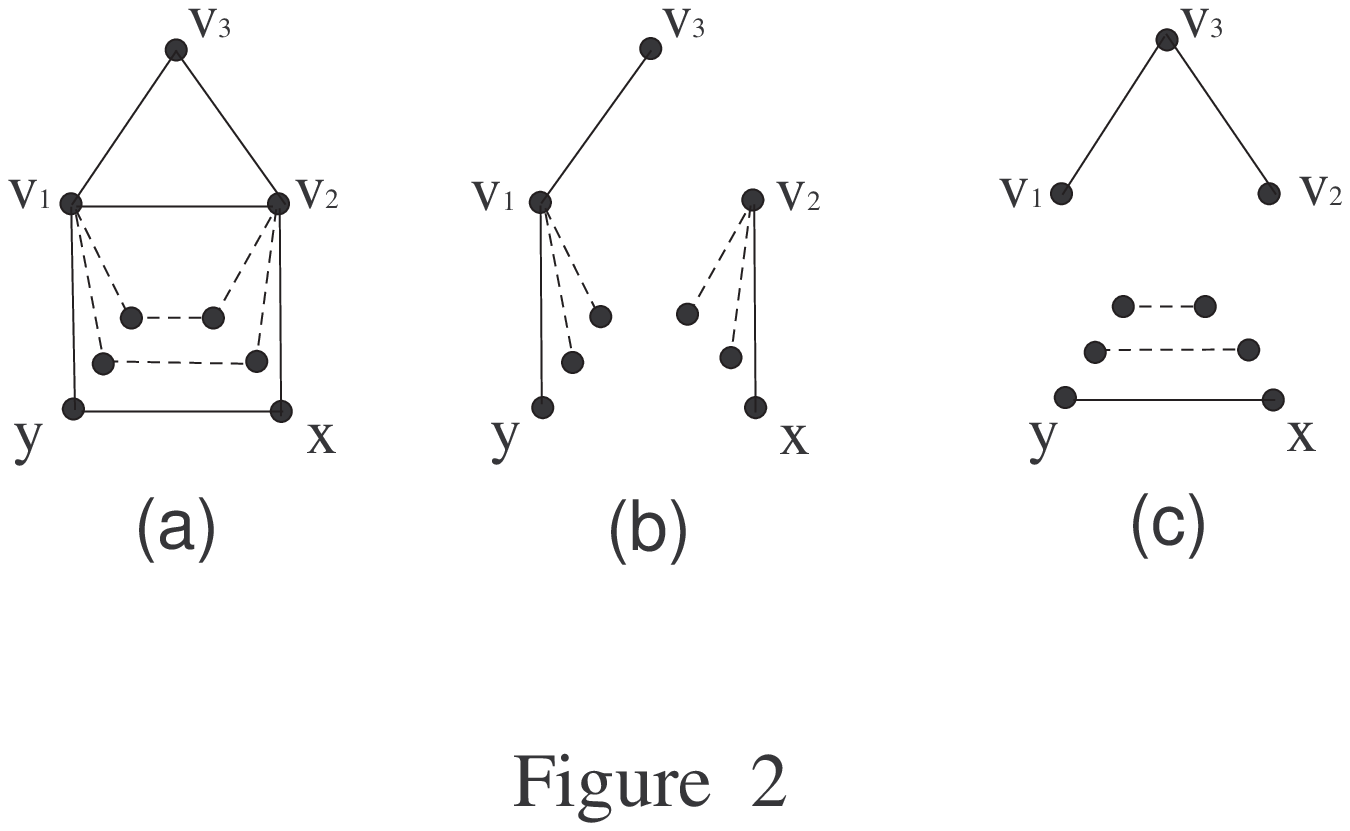}
\end{center}
\end{figure}

If $d_{G}(v_{1})=d_{G}(v_{2})=3$, then $G\in \mathscr{U}$.
Otherwise $d_{G}(v_{1})=d_{G}(v_{2})=k \ (k\geq 4)$, then $G$
contains two star-factors shown in Figure $2(b)$ and $2(c)$, with
the weights $2k-3$ and $k$, respectively. However, we see
$2k-3\neq k$ for $k\geq 4$, a contradiction to $G\in \mathscr{U}$.
Hence there is exactly one quadrangle using the edge $v_{1}v_{2}$
in $G$, and $G$ is shown in Figure $1(b)$.

\vspace{2mm}

{\it Case $2$.} $d_{G}(v_{i})\geq 3, i=1,2,3$.

\vspace{2mm}

{\it Claim 1.} Let $X=\{x\in V(G) \mid N_{G}(x)\subseteq
\{v_{1},v_{2},v_{3}\}\}$, then $X=\emptyset$.

\vspace{2mm}

Suppose $X\neq \emptyset$. Let $F_{3}=\{u_1u_2\ | \ u_1\in \{
v_{1},v_{2},v_{3}\}, u_2\in V(G)-X-\{ v_{1},v_{2},v_{3}\}\}$. Then
$G-F_{3}$ has no isolated vertices, so $G-F_{3}\in \mathscr{U}$ by
Lemma $1$. Assume

$$N_{G}(v_{1})\cap N_{G}(v_{2})=\{x_{1}, \cdots , x_{i}\},$$
$$N_{G}(v_{2})\cap N_{G}(v_{3})=\{y_{1}, \cdots , y_{j}\},$$
$$N_{G}(v_{1})\cap N_{G}(v_{3})=\{z_{1}, \cdots , z_{k}\},$$
$$N_{G}(v_{1})\cap N_{G}(v_{2})\cap N_{G}(v_{3})=\{u_{1}, \cdots , u_{l}\}.$$

Then $G-F_{3}$ contains a component $H$ with vertices in $X$ and the
triangle $\triangle v_{1}v_{2}v_{3}$ shown in Figure $3(a)$. Since
$X\neq \emptyset$, without loss of generality, we assume $l=0$ and
at least one of $i,j,k$ is nonzero.

{\it Subcase 1.1.} There is exactly one of $i,j,k$ is nonzero. Then
$H$ can be decomposed into one or two stars.

{\it Subcase 1.2.} There are exactly two of $i,j,k$ are nonzero.
Assume, without loss of generality, that $i=0$, $j\neq 0$ and $k\neq
0$. If $j = k =1$,  then $H$ can be decomposed into one or two
stars. Otherwise $H$ can be decomposed into one or two or three
stars.

{\it Subcase 1.3.} $i,j,k$ are all nonzero. Then $H$ can be
decomposed into two or three stars.

 So in all three subcases it contradicts to $G-F_{3}\in
\mathscr{U}$.

\vspace{3mm}

\begin{figure}[h,t]
\begin{center}
\includegraphics[width=11cm]{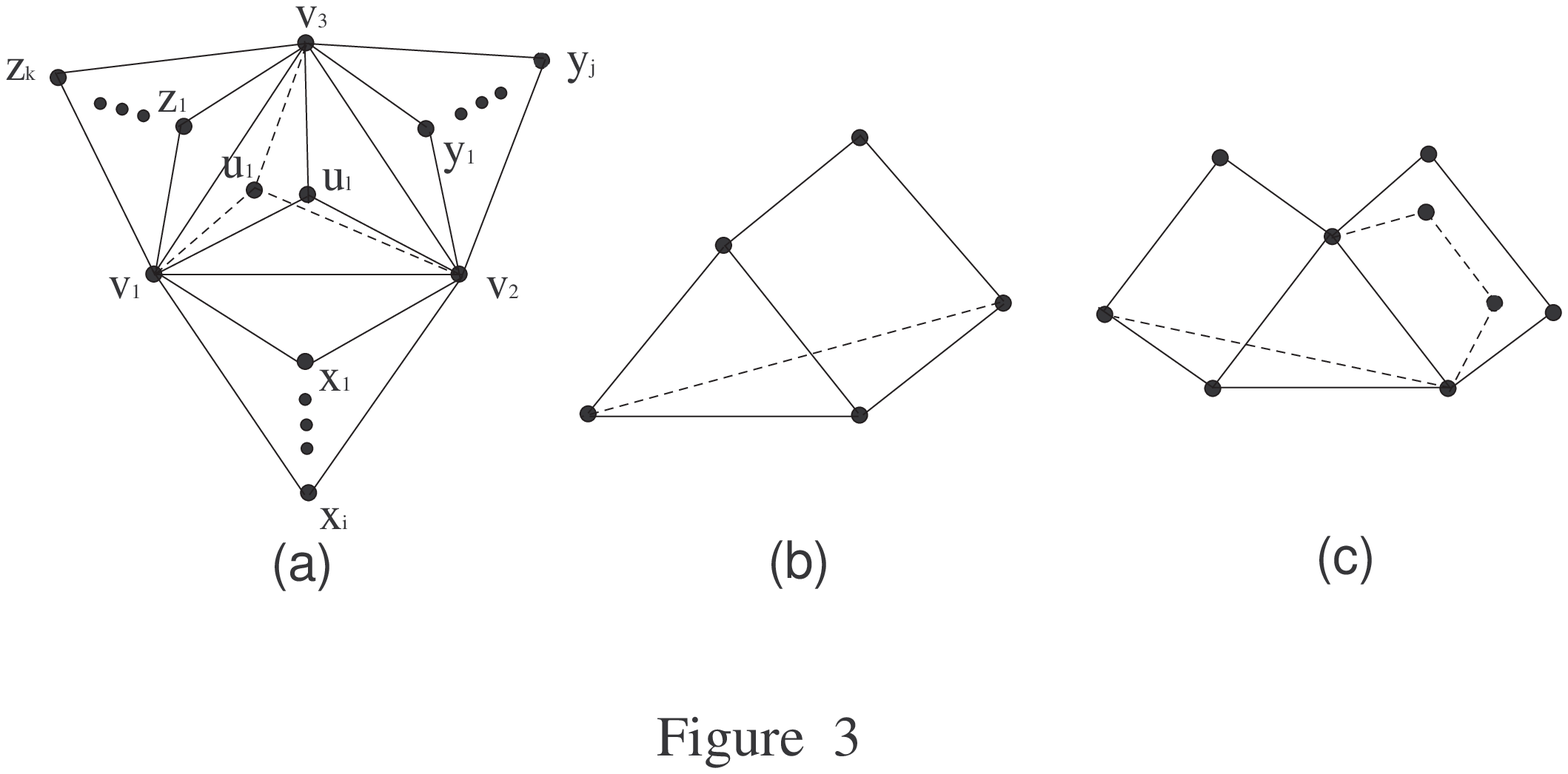}
\end{center}
\end{figure}

Let $F_{4}=\{u_1u_2 \ | \ u_1\in \{v_{1},v_{2}\}, u_2\in
V(G)-\{v_{1},v_{2},v_{3}\}\}$. Then, by Claim $1$, no isolated
vertices are created in $G-F_{4}$, and all the neighbors of $v_{3}$
other than $v_{1}$ and $v_{2}$ are stems in graph $G-F_{4}$ by Lemma
$3$. Suppose $u$ is a stem which is adjacent to $v_{3}$ in
$G-F_{4}$, and $m$ is a leaf adjacent to $u$ in $G-F_{4}$. It is
obvious that $m$ can only be adjacent to $v_{1}$ or $v_{2}$ besides
$u$ in $G$.

\vspace{2mm}

{\it Claim 2.} $N_{G}(u)-\{m, v_{3}\}\subseteq \{v_{1}, v_{2}\}$
and $N_{G}(u)\cap N_{G}(m)=\emptyset$.

\vspace{2mm}

{\it Subcase 2.1.} $m$ is adjacent to exactly one of $v_{1}$ and
$v_{2}$ in $G$. Assume, without loss of generality, that $m$ is
adjacent to $v_{1}$, and $u$ has other neighbors other than $v_{1},
v_{2}, v_{3}$ and $m$. Let $F_{5}=\{u_1u_2 \ | \ u_1\in
\{v_{2},v_{3}\}, u_2\in V(G)-\{v_{1},v_{2},v_{3}\}\}$, then $u$ is
not a leaf in $G-F_{5}$. But the only neighbor of $m$ in $G-F_{5}$
is $u$, so $m$ is not a stem in $G-F_{5}$. However, $m$ should be a
stem in $G-F_{5}$ by Lemma $3$, a contradiction. Hence
$N_{G}(u)-\{m, v_{3}\}\subseteq \{v_{1}, v_{2}\}$. If $u$ is also
adjacent to $v_{1}$, then $m$ and $u$ are stems in graph $G-F_{5}$
by Lemma $3$, a contradiction.

{\it Subcase 2.2.} $m$ is adjacent to both $v_{1}$ and $v_{2}$ in
$G$. One may obtain a contradiction by a similar argument as in
Subcase 2.1.\\

    Hence $v_{1}$, $v_{2}$, $v_{3}$, $m$ and $u$ form at most two
quadrangles with common edge $mu$ by Claim $2$, and an induced
subgraph of $G$ with the vertices $v_{1}$, $v_{2}$, $v_{3}$, $m$ and
$u$ is isomorphic to the graph shown in Figure $3(b)$ (dashed line
indicates a possible edge). So the subgraph $H$ induced by all the
vertices in the component, in $G-F_{4}$, which contains the
$3$-cycle $C_{3}=v_{1}v_{2}v_{3}$ is isomorphic to the graph shown
in Figure $3(c)$.

    Let $F_{6}=\{u_1u_2 \ | \ u_1\in \{v_{1},v_{3}\}, u_2\in
V(G)-\{v_{1},v_{2},v_{3}\}\}$. By the similar argument above, both
subgraphs induced by the vertices in the component of $G-F_{5}$ and
$G-F_{6}$, respectively, which contain the $3$-cycle
$C_{3}=v_{1}v_{2}v_{3}$ are also isomorphic to the graph shown in
Figure $3(c)$. So $G$ is isomorphic to the graph shown in Figure
$4(a)$. If we delete some edges from $G$ such that all vertices in
$G-\{v_{1}, v_{2}, v_{3}\}$ are of degree two, then the spanning
subgraph $G'$ of $G$ will be the graph shown in Figure $4(b)$.

\vspace{3mm}

\begin{figure}[h,t]
\begin{center}
\includegraphics[width=9cm]{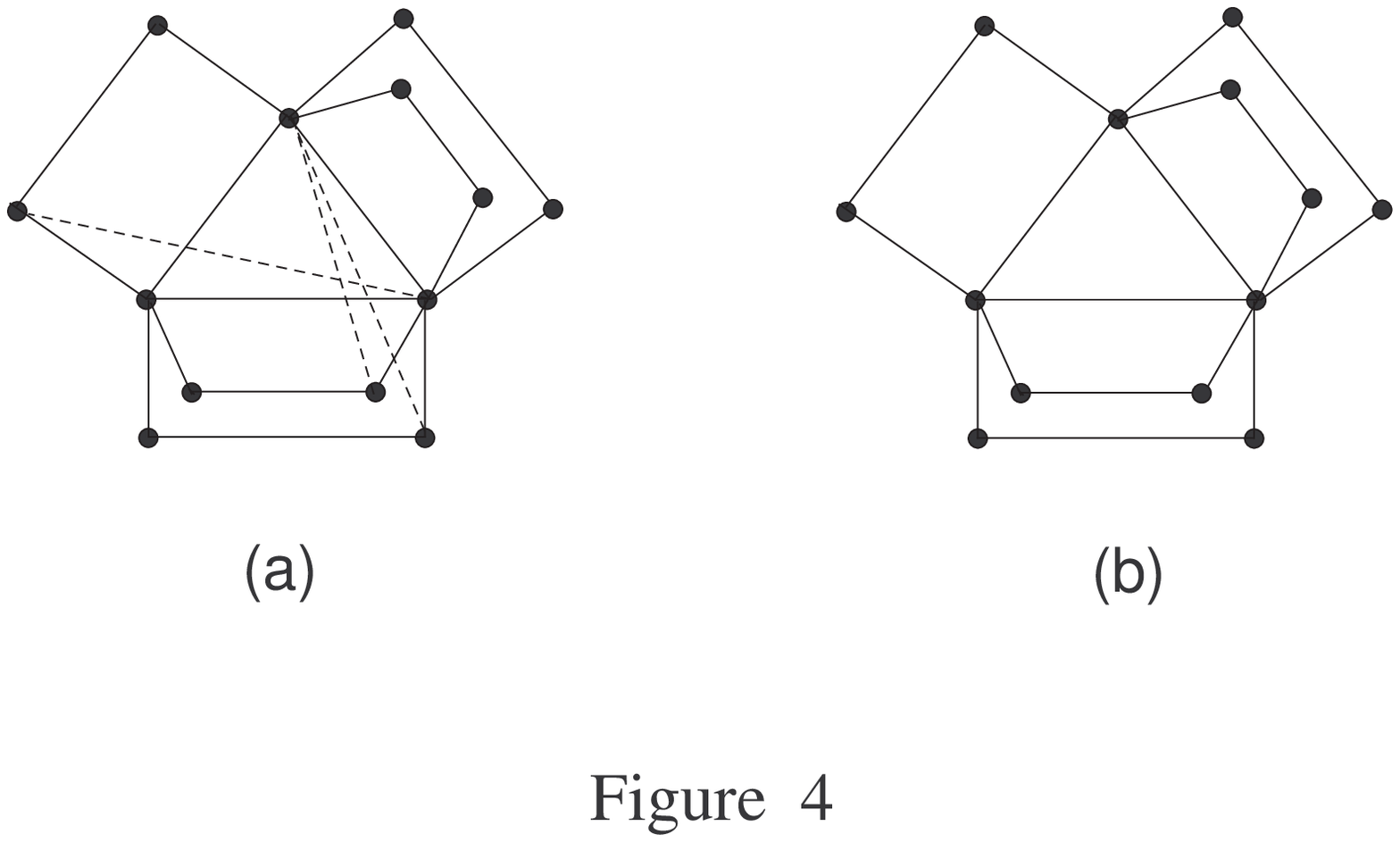}
\end{center}
\end{figure}

{\it Claim 3.} For each edge on the triangle $\Delta
v_{1}v_{2}v_{3}$, there exits at most one quadrangle in $G'$
containing it.

\vspace{2mm}

Let $Y$ denote the vertices that are contained in a quadrangle which
use the edge $v_{1}v_{2}$ and $F_{7}=\{u_1u_2 \ | \ u_1\in
\{v_{1},v_{2},v_{3}\}, u_2\in V(G)-Y-\{v_{3}\}\}$. Then a component
containing the triangle $\Delta v_{1}v_{2}v_{3}$ of $G-F_{7}$ is
either the triangle $\Delta v_{1}v_{2}v_{3}$ itself or a triangle
satisfying the conditions in Case $1$. So $v_{1}v_{2}$ is contained
in at most one quadrangle in $G$. The same argument can be applied
to edges $v_{2}v_{3}$ and $v_{1}v_{3}$.

\vspace{3mm}

\begin{figure}[h,t]
\begin{center}
\includegraphics[width=10cm]{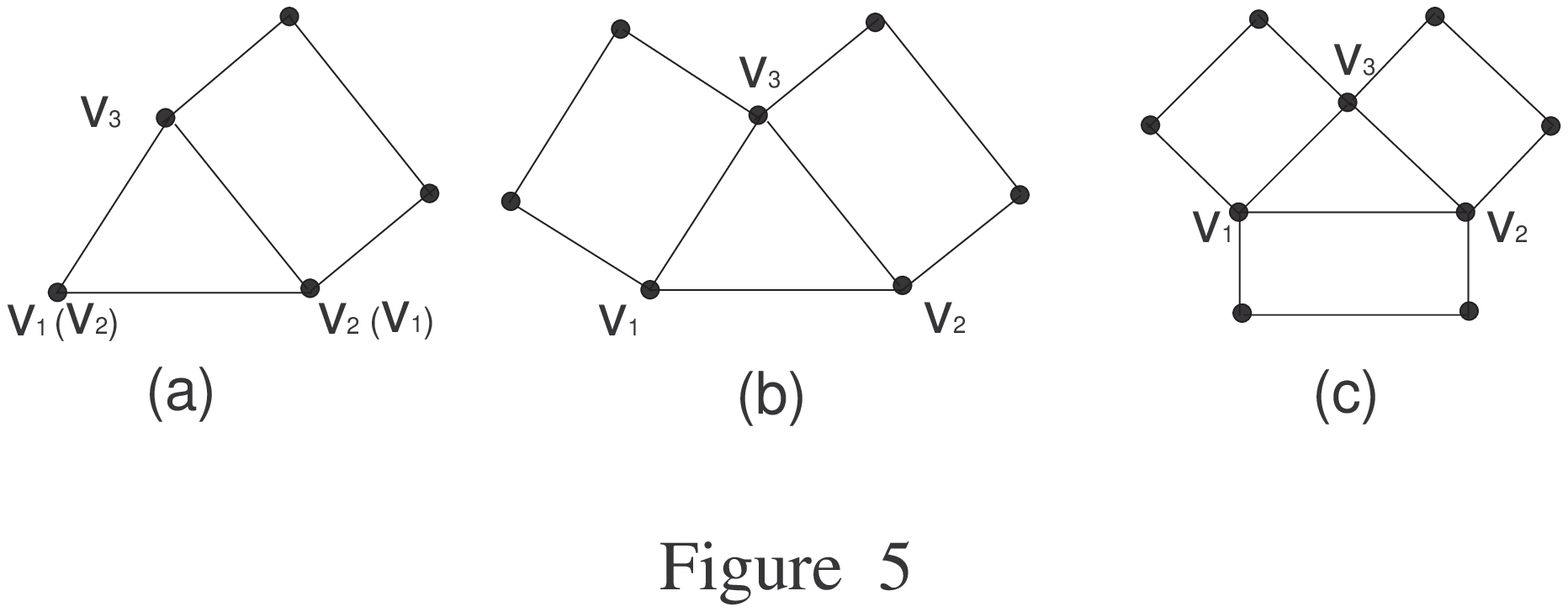}
\end{center}
\end{figure}

    Now we know that $G'$ could be the graphs shown in Figure $5$.
However, the graph shown in Figure $5(c)$ can be decomposed into
three or four stars. Hence, by Lemma 1, $G \notin \mathscr{U}$. So
the only possible graph $G'$ are the graphs shown in Figure $5(a)$
and $5(b)$. We add the edges back following the principle of Claim
$2$, then $G$ can only be the graphs shown in Figure $1(c)$, $1(d)$
and $1(e)$ since every vertex in the triangle $\Delta
v_{1}v_{2}v_{3}$ has degree at least three in $G$.

    This completes the proof of Theorem 2.\hfill$\Box$\\

    The main theorem has classified all graphs in $\mathscr{U}$ with
girth three. Combining with Theorem 1, the only two families
remaining to be determined are graphs of girth four or graphs with
leaves and small girths. It seems that the structures of both
families are much more complicated and new techniques are required
in order to determine them completely.


\begin{thebibliography}{99}

\bibitem{amahashi} A. Amahashi and M. Kano,
On factors with given components, {\em Discrete Math.}, 42(1982),
1-6.

\bibitem{bb} B. Bollob\'{a}s, {\it Modern Graph Theory}, 2nd Edition,
Springer-Verlag New York, Inc. 1998.

\bibitem{hartnell} B. L. Hartnell and D. F. Rall,
On graphs having uniform size star factors, {\em Australas. J.
Combin.}, 34(2006), 305-311.

\bibitem{yu} Q. Yu,
Counting the number of star-factors in graphs, {\em J. Combin.
Math. Combin. Comput.}, 23(1997), 65-76.


\end{thebibliography}
\end{document}